\documentclass[12pt]{article}
\usepackage{latexsym}
\usepackage{amssymb}
\usepackage[T1]{fontenc}
\fontencoding{T1}
\usepackage[latin2]{inputenc}

\def\vsp{\vspace{2mm}}

\def\Nn{\mathbb{N}}
\def\Pp{{\bf P}}
\def\Rr{\mathbb{R}}
\def\Ee{{\bf E}}

\def\qed{\hfill{$\Box $}}

\newtheorem{cor}{\large\bf Corollary}
\newtheorem{prop}{\large\bf Proposition}
\newtheorem{definition}{\large\bf Definition}
\newtheorem{theorem}{\large\bf Theorem}
\newtheorem{lemma}{\large\bf Lemma}

\newenvironment{namelist}[1]{%
\begin{list}{}
     {
      
      \settowidth{\labelwidth}{#1}
      \setlength{\leftmargin}{1.1\labelwidth}
               }
      }{%
\end{list}}

\begin{document}
\title{On weak generalized stability \newline and $(c,d)$-pseudostable
random variables \newline via functional equations \newline {\small Dedicated to the memory of Andrzej Lasota}}
\author{W. Jarczyk and J. K. Misiewicz
\thanks{Department of Mathematics Informatics and Econometrics,
University of Zielona G{\'o}ra,
 ul. Szafrana 4A, 65-516 Zielona G{\'o}ra,
 Poland \newline
This paper was partially written while the second author was a visiting professor of Delft Institute of Applied Mathematics, Delft University of Technology, Holland}}
\date{}
\maketitle

\begin{abstract} In this paper we give a first attempt to define and study stable distributions with respect to
the weak generalized convolution, focusing our attention on the symmetric weakly stable distribution. As in the
case of the classical convolution, characterization of distributions stable in the sense of the weak generalized
convolution depends on solving some functional equations in the class of characteristic functions.
\end{abstract}

\noindent
{\bf Key words}: weakly stable distribution,
symmetric stable distribution, \linebreak
$\ell_\alpha$-symmetric distribution,
generalized convolution, scale mixture, functional equation,
convexity \\
{\bf Mathematics Subject Classification:} 60A10, 60B05, 60E05,
60E07, 60E10, 39B22, 26A51.

\section{Introduction}

The investigations of weakly stable random variables started in the seventies in the papers of Kucharczak and
Urbanik (see \cite{KU, Urbanik1}).  Later  a series of papers on weakly stable distributions written by Urbanik,
Kucharczak and Vol'kovich was published (see e.g. \cite{KU2, vol1, vol2, vol3}). Recently  a paper written by
Misiewicz, Oleszkiewicz and Urbanik (see \cite{MOU}) appeared, where one can find a full characterization of
weakly stable distributions with non-trivial discrete part, and a substantial attempt to characterize weakly
stable distributions in the general case.

In stochastic modeling of real processes, using independent random variables or Gaussian processes in a variety
of constructions turned out to be not sufficient or adequate. Multidimensional stable distributions have nice
linear properties and enable more complicated structures of dependencies, thus recently the role of stable
processes in stochastic modeling is growing. On the other hand stable distributions are very difficult in
calculations. However there are some efficient techniques for their computer simulation.

In this situation weakly stable distributions and processes seem to be good candidates for use in stochastic
modeling. They extend Feller's idea of subordinated processes. They have nice linear properties, namely: if
$({\bf X}_i)$ is a sequence of independent identically distributed random vectors with a weakly stable
distribution, then every linear combination $\sum a_i {\bf X}_i$ has the same distribution as ${\bf X}_1 \cdot
\theta$ for some random variable $\theta$ independent of ${\bf X}_1$. This condition holds not only when $(a_i)$
is a sequence of real numbers, but also when $(a_i)$ is a sequence of random variables such that $(a_i)$ and
$({\bf X}_i)$ are independent. This means that dependence structure of the linear combination $\sum a_i {\bf
X}_i$ and dependence structure of the random vector ${\bf X}_1$ are the same, and the sequence $(a_i)$ is
responsible only for the radial behavior. Moreover, weak stability is preserved under taking linear operators,
projections or functionals. On the other hand, radial properties of distribution can be arbitrarily defined by
choosing a proper random variable $\theta$ independent of ${\bf X}_1$ and considering the distribution of
$\theta \cdot {\bf X}_1$. Similar properties of tempered stable distributions (see e.g. \cite{JR}) are the
reason why they are so important now in statistical physics  modelling turbulence, or in mathematical finance for
modelling  stochastic volatility.

In this paper we develop the idea of distributions stable with respect to a generalized convolution defined by
weakly stable variables. In this construction the weakly stable variable plays the role of a catalyst in the
presence of which the underlining process can develop. This weakly stable variable can be also treated as a
filter, so that we observe the original process only by its filtered values. In this sense, we want to
characterize distributions which are stable after filtering.

By $\mathcal{P}(\mathbb{E})$ we denote the set of all probability measures on a separable Banach space
$\mathbb{E}$ (with dual $\mathbb{E}^{\ast}$).  For simplicity we write $\mathcal{P}$ for the set of all
probability measures on $\Rr$. Moreover, $\mathcal{P}_{+}$ denotes the set of probability measures on
$[0,\infty)$. The symbol $\Phi$ stands for the set of all characteristic functions on $\mathbb{R}$, whereas
$\delta_{\mathbf{x}}$  denotes the probability measure concentrated at the point $\mathbf{x} \in \mathbb{E}$. If
a sequence of probability measures $(\lambda_n)_{n\in \Nn}$ converges weakly to a probability measure $\lambda$,
we write $\lambda_n \rightarrow \lambda$.

For a random vector (or random variable) $\mathbf{X}$, we write $\mathcal{L}(\mathbf{X})$ for  the
distribution of $\mathbf{X}$. For $\lambda \in \mathcal{P}$, $\lambda = \mathcal{L}(\theta)$, we write
$|\lambda| = \mathcal{L}( |\theta|)$. If $\mu = \mathcal{L}(\mathbf{X})\in \mathcal{P}(\mathbb{E})$, then the
characteristic function $\widehat{\mu} \colon \mathbb{E}^{\ast} \rightarrow \mathbb{C}$ of the measure $\mu$ (of
the random vector $\mathbf{X}$) is defined by
$$
\widehat{\mu}(\xi) = \Ee \exp \left\{ i <\xi, \mathbf{X}> \right\} = \int_{\mathbb{E}} \exp\left\{ i <\xi, x >
\right\} \mu(dx).
$$
 For random vectors $\mathbf{X}, \mathbf{Y}$, we
 write $\mathbf{X} \stackrel{d}{=} \mathbf{Y}$ for $\mathcal{L}(\mathbf{X}) = \mathcal{L}
(\mathbf{Y})$. If $\mathbf{X}$ and $\mathbf{Y}$ are independent random vectors then $\mathcal{L}(\mathbf{X} +
\mathbf{Y})$ is the convolution of $\mathcal{L}(\mathbf{X})$ and $\mathcal{L}(\mathbf{Y})$,  denoted by
$\mathcal{L}(\mathbf{X}) \ast \mathcal{L}(\mathbf{Y})$.

For $t \in \Rr$, a rescaling operator $T_t \colon \mathcal{P}(\mathbb{E}) \rightarrow \mathcal{P}(\mathbb{E})$ is
defined as follows:
$$
T_t \mu (A) = \left\{ \begin{array}{lcl}
            \mu({A/t}) & \hbox{ if } & t \in \Rr\setminus
              \{0\}, \\
            \delta_0(A) & \hbox{ if } & t=0.
            \end{array} \right.
$$
It is easy to see that if $\mu = \mathcal{L}(\mathbf{X})$ then
$T_t \mu = \mathcal{L}(\mathbf{tX})$. The scale mixture $\mu
\circ \lambda$ of the measure $\mu \in \mathcal{P}(\mathbb{E})$
with respect to the measure $\lambda \in \mathcal{P}$ is
defined by
$$
\mu \circ \lambda (A) = \int_{\Rr} T_t \mu(A) \lambda(dt).
$$
If $\mu = \mathcal{L}(\mathbf{X})$ and $\lambda = \mathcal{L}(\theta)$ with $\mathbf{X}$ and $\theta$
independent, then $\mu \circ \lambda = \mathcal{L} (\mathbf{X} \theta)$.

\vspace{3mm}

A random vector ${\bf X}$ with the distribution $\mu$ on a real
separable Banach space $\mathbb{E}$ is weakly stable iff
$$
\forall \, a,b \in \Rr \,\, \exists \, \theta \hspace{5mm} a{\bf X} + b {\bf X}' \stackrel{d}{=} {\bf X} \theta,
\eqno{(A)}
$$
where ${\bf X}'$ is an independent copy of ${\bf X}$ and the random variable $\theta$ is independent of ${\bf X}$.
It was proved in \cite{MOU} that ${\bf X}$ is weakly stable if and only if
$$
\forall \, \theta_1, \theta_2 \,\, \exists \, \theta \hspace{5mm} \theta_1 {\bf X} + \theta_2 {\bf X}'
\stackrel{d}{=} {\bf X} \theta, \eqno{(B)}
$$
where $\theta_1, \theta_2$ are real random variables such that $\theta_1, \theta_2, {\bf X}, {\bf X}'$ are
independent and the random variable $\theta$ is independent of ${\bf X}$. In the language of probability
measures, the condition $(B)$ can be written as
$$
\forall \, \lambda_1, \lambda_2 \in \mathcal{P} \,\, \exists \, \lambda \in \mathcal{P} \hspace{5mm} (\mu \circ
\lambda_1) \ast (\mu \circ \lambda_1) = \mu \circ \lambda, \eqno{(C)}
$$
for $\lambda, \lambda_i$ the distributions of $\theta, \theta_i$, $i=1,2$. It was shown in \cite{MOU} that
the measure $\lambda$ is uniquely determined if the measure $\mu$ is not symmetric. For a symmetric measure $\mu$,
we only have uniqueness of the measure $|\lambda|$.

The best known examples of weakly stable random vectors are symmetric stable vectors, and in this case the random
variable $\theta$ appearing in the condition $(A)$ is a constant, $\theta \equiv (|a|^{\alpha} +
|b|^{\alpha})^{1/{\alpha}}$ for some $\alpha \in (0,2]$. Another family of weakly stable distributions consists
of uniform distributions $\omega_n$ on unit spheres $S_{n-1} \subset \Rr^n$ and their lower-dimensional
projections.

It was shown in \cite{MOU} that if a weakly stable distribution $\mu$ contains a discrete part, then it is
discrete and either $\mu = \delta_0$, or $\mu = \frac{1}{2} \delta_a + \frac{1}{2} \delta_{-a}$ for some $a \in
\mathbb{E}\smallsetminus \{0\}$. From now on we will assume that the considered weakly stable measure $\mu$ is
non-trivial in the sense that it is not discrete.

We can now define a generalized weak convolution $\oplus = \oplus_{\mu}$ for any nontrivial weakly stable
measure $\mu$.

\begin{definition} Let $X$ be a non-trivial random vector with the weakly stable distribution $\mu$.  The weak
generalized convolution $\oplus$ of measures $\lambda_1,\lambda_2 \in \mathcal{P}$  is defined by
$$
\lambda_1 \oplus \lambda_2 = \left\{ \begin{array}{ll}
   \lambda & \hbox{ if $\mu$ is not symmetric}, \\
   |\lambda| & \hbox{ if $\mu$ is symmetric},
   \end{array} \right.
$$
where $\lambda \in \mathcal{P}$ is such that $(\mu \circ \lambda_1) \ast (\mu \circ \lambda_2) = \mu \circ
\lambda$. For two independent random variables $\theta_1$ and $\theta_2$ with distributions $\lambda_1$ and
$\lambda_2$ respectively, the weak generalized sum $\theta_1 \oplus \theta_2$ is the random variable defined by
$$
\theta_1 {\bf X} + \theta_2 {\bf X}' \stackrel{d}{=} {\bf X} \left( \theta_1 \oplus \theta_2 \right),
$$
where $\theta_1, \theta_2, {\bf X}, {\bf X}'$ are independent, the random variable $\theta_1 \oplus \theta_2$ is
independent of ${\bf X}$ and $\mathcal{L}(\theta_1 \oplus \theta_2) = \lambda_1 \oplus \lambda_2$.
\end{definition}

The operation $\oplus$ in $\mathcal{P}$ is commutative and associative. Moreover, as shown in
\cite{Mis06}, the following conditions hold:
\begin{namelist}{ll}
\item[(i)] the measure $\delta_0$ is the unit element, i.e. $\delta_0 \oplus \lambda = \lambda$ for all $\lambda
     \in {\cal P}$ if $\mu$ is not symmetric and $\delta_0 \oplus \lambda = |\lambda |$ for all $\lambda \in
     {\cal P}$ if $\mu$ is symmetric;
\item[(ii)] $(p \lambda_1 + q \lambda_2) \oplus \lambda = p
     (\lambda_1 \oplus \lambda) + q (\lambda_2 \oplus
     \lambda)$, whenever $\lambda, \lambda_1, \lambda_2 \in
     {\cal P}$ and $p\geqslant 0$, $q\geqslant 0$, $p+q = 1$
     (linearity);
\item[(iii)] $(T_a \lambda_1) \oplus (T_a \lambda_2) = T_a
      (\lambda_1 \oplus \lambda_2)$ for any $\lambda_1,
      \lambda_2 \in {\cal P}$ and $a>0$ (homogeneity);
\item[(iv)] if $\lambda_n \rightarrow \lambda_0$ then $\lambda_n
     \oplus \lambda \rightarrow \lambda_0 \oplus \lambda$
     for all $\lambda \in {\cal P}$ (continuity).
\end{namelist}

The idea of generalized convolutions has been extensively studied after it was introduced by K. Urbanik in 1964
\cite{Urbanik64}. The definition proposed by K. Urbanik is as follows:

A commutative and  associative binary operation $\diamond : {\cal P}_{+} \times {\cal P}_{+} \rightarrow {\cal
P}_{+}$ is called a {\em generalized convolution} if it satisfies conditions (i)$\div$(iv) with $\oplus$
replaced by $\diamond$ and the following condition holds:
\begin{namelist}{ll}
\item[(v)] there exists a sequence $(c_n)_{n\in \Nn}$ of positive
     numbers such that the sequence $(T_{c_n}
     \delta_1^{\diamond n})_{n \in \Nn}$ weakly converges to a
     measure different from $\delta_0$.
\end{namelist}

\vspace{2mm}

The first, but not the most important difference between the definition of generalized convolution $\diamond$
given by K.  Urbanik and the definition of weak generalized convolution $\oplus$ is the domain, i.e. $\oplus :
{\cal P} \rightarrow {\cal P}$. This implies, in particular, that most of the methods used in studying
generalized convolutions cannot be directly applied for weak generalized convolutions. Another difference is
that the weak generalized convolution have properties (i), (ii), (iii) and (iv), but, in general, not (v). In
spite of this disadvantage, we do not have to assume that the algebra $({\cal P}, \oplus_{\mu})$ is regular,
i.e. that there exists a non-trivial homomorphism of $({\cal P}, \oplus_{\mu})$ into a complex field. The
assumption that there exists a non-trivial homomorphism from $({\cal P}_{+}, \diamond)$ into the positive
half-line was crucial in studying generalized convolutions. In the case of the weak generalized convolution
$\oplus_{\mu}$ we have that for every $\xi \in \mathbb{E}^{\ast}$ the formula
$$
h_{\xi}(\lambda) = \int_{\mathbb{R}} \widehat{\mu}(t \xi ) \lambda(dt)
$$
defines a homomorphism of $({\cal P}, \oplus_{\mu})$ into a
complex field. Moreover, if the weakly stable measure $\mu$ is
non-trivial, then there exists $\xi \in \mathbb{E}^{\ast}$ such that $h_{\xi}$ is non-trivial.

Talking about the homomorphism and treating the set $\mathcal{P}$ as an algebra we underline here that
$\mathcal{P}$ is equipped with the generalized convolution $\oplus_{\mu}$ as a binary operation and with
rescaling measures operator $T_t $, $t\in \mathbb{R}$, which can be treated as multiplication by scalars.
Moreover, we see that convex linear combinations of probability measures are also probability measures, thus
$\mathcal{P}$ can be treated as a subset of some linear space.

In this paper we consider only symmetric weakly stable distributions $\mu$, thus the functions $h_{\xi}$, $\xi
\in \mathbb{E}^{\ast}$, and $\widehat{\mu}$ are real-valued. Moreover, for symmetric $\mu$ we know that $\mu
\circ \lambda = \mu \circ |\lambda|$, so we can restrict our attention to the set $(\mathcal{P}_{+},
\oplus_{\mu})$
instead of $(\mathcal{P}, \oplus_{\mu})$. \\

\vspace{0.1mm}

 The paper is organized as follows. In section 2 solve a functional equation, which will be needed
later. In section 3 we describe the problem of characterizing stable distributions with respect to the
generalized convolution $\oplus_{\mu}$. The description of strictly stable distributions in this sense is given
in subsection 3.1. In the following subsections we discuss the possibility of having occasionally some strange
behavior of the function $d$. The last subsection contains discussion of the general case of stable, but not
strictly stable, distributions in the sense of the generalized convolution $\oplus_{\mu}$.

\section{Functional equation}

Let $\mathcal{F}$ be the set of all continuous functions $f \colon [0,\infty) \rightarrow \Rr$ such that $f(0) =
0$ and $0$ is an isolated point of the set $\{ x\in \Rr\colon f(x) = 0\}$.

\begin{theorem}
Let $a,b \in (0,\infty)$ and let $f \in \mathcal{F}$ be a solution of the functional equation
$$
 f(t) = f(at) + f(bt). \eqno{(1)}
$$
Then $a,b \in (0,1)$ and there exist $p > 0$ and a continuous function $H \colon (0,\infty) \rightarrow
\Rr\setminus \{0\}$ such that $H(t) = H(at) = H(bt)$ for every $t > 0$, and
$$
 f(t) = t^{p} H(t),  \hspace{5mm} t>0.
$$
Moreover, if ${{\ln a}/{\ln b}}$ is irrational then the
function $H$ is constant.
\end{theorem}

\vspace{2mm}

The proof of Theorem 1 is based on a series of lemmas. The first one can be easily proved by mathematical
induction.

\vspace{2mm}

\noindent
\begin{lemma}
Let $a,b \in (0,\infty)$ and let $f \colon [0,\infty) \rightarrow \Rr$ be a solution of equation $(1)$. Then for
every $n \in \Nn$ and every $t\in \Rr$ we have
$$
f(t) = \sum_{k=0}^n {{n}\choose{k}} f\left( a^k b^{n-k} t \right).
$$
\end{lemma} \qed

\vspace{2mm}

\noindent
\begin{lemma}
Let $a,b \in (0,\infty)$. If $(1)$ has a solution in the class $\mathcal{F}$ then $a,b \in (0,1)$.
\end{lemma}

\vspace{2mm}

\noindent {\bf Proof.} Let $f \in \mathcal{F}$ be a solution of $(1)$. Without loss of generality we may assume
that $b\leqslant a$ and $f(x) >0$ for every $x \in (0,u)$ with some $u>0$. Then for every $t\in (0,u)$ we have
$ba^{-1}t \in (0,u)$ and, consequently,
$$
 f(t) = f(a a^{-1} t) = f(a^{-1}t) - f(ba^{-1}t) < f(a^{-1} t).
$$
If $a>1$ this and mathematical induction would imply that $f(t) \leqslant f(a^{-n}t)$ for every $n\in
\mathbb{N}$,
 and by the continuity of $f$ we would have $f(t) = 0$ for every $t\in (0,u)$. If $a=1$ then
equation $(1)$ imply that $f(bt) = 0$ for each $t>0$. In both cases $f \not\in \mathcal{F}$ which contradicts
our assumption. Therefore $a<1$. \qed

\vspace{2mm}

\begin{lemma}
Let $a,b \in (0,\infty)$ and let $f \in \mathcal{F}$ be a solution of equation $(1)$. Then $f(t) = 0$ if and
only if $t=0$.
\end{lemma}

\vspace{2mm}

\noindent {\bf Proof.} Assume, for instance, that for some $u>0$ we have $f(x) > 0$ for every $x \in (0,u)$. By
Lemma 2 we know that $a,b \in (0,1)$, thus for every fixed $t>0$ we can find $n \in \Nn$ large enough to have
$$
a^k b^{n-k} t \in (0,u), \hspace{5mm}  k = 0,1,\dots, n.
$$
Then, by Lemma 1, we have
$$
f(t) = \sum_{k=0}^n {{n}\choose{k}} f\left( a^k b^{n-k} t
\right) > 0
$$
which ends the proof. \qed

\vspace{2mm}

 The proofs of Lemmas 4--6 below take pattern of some ideas
from \cite{Jar} (see also \cite{BJ}).

\vspace{2mm}

\begin{lemma}
Let $i,j \in \mathbb{Z}$. For every $\alpha >0$ and $\beta \in (0,2\alpha)$ there exists $n_0 \geqslant \max \{
|i|, |j|\}$ such that for every $n\geqslant 2n_0$
$$
{{n+i+j}\choose{k+i}} + \alpha^2  {{n-i-j}\choose{k-i}} \geqslant \beta {n \choose k}, \hspace{8mm} k \in \{n_0,
\dots, n-n_0\}.
$$
\end{lemma}

\vspace{2mm}

\noindent {\bf Proof.} Let $\alpha >0$ and $\beta \in (0, 2\alpha)$ be fixed. Assume first that $i+j \geqslant
0$. Then at least one of the numbers $i,j$, for example $i$, is nonnegative. There exists $n_0 \in \Nn$ such
that $n_0 \geqslant \max\{ |i|, |j|\}$ and
$$
\left( 1 - \frac{i}{n_0 +1} \right)^i \left( 1- \frac{|j|}{n_0 + 1}\right)^{|j|} \geqslant \frac{\beta}{2
\alpha}.
$$
Fix integers $n$ and $k$ such that $n \geqslant 2n_0$ and $n_0 \leqslant k \leqslant n-n_0$. Then we have
\begin{eqnarray*}
\lefteqn{ {{n+i+j}\choose{k+i}} + \alpha^2
       {{n-i-j}\choose{k-i}} = {n \choose k} \left[ w_n + \frac{
       \alpha^2}{w_n}\cdot \frac{ (n+1) \dots (n+i+j)}{ (n-i-j +1)
       \dots n} \right. } \\
&& \left. \cdot \left(1-\frac{i}{k+1} \right) \dots \left(1 -
\frac{i}{k+i} \right) \left( 1- \frac{|j|}{n-k+1} \right) \dots
\left( 1- \frac{|j|}{n-k +|j|} \right) \right],
\end{eqnarray*}
where
$$
w_n = \frac{(n+1) \dots (n+i+j) \cdot (n-k)!}{ (k+1) \dots
(k+i) \cdot (n-k+j)!}.
$$
This implies that
\begin{eqnarray*}
\lefteqn{ {{n+i+j}\choose{k+i}} + \alpha^2
       {{n-i-j}\choose{k-i}} } \\
&& \geqslant {n \choose k} \left[ w_n +
       \frac{\alpha^2}{w_n} \left( 1 - \frac{i}{n_0 +1}
       \right)^i \left( 1- \frac{|j|}{n_0 + 1}\right)^{|j|}
       \right] \\
&& \geqslant {n \choose k} \left[ w_n +
       \frac{\beta^2}{4} \frac{1}{w_n} \right]  \geqslant {n \choose
       k} \left[ \frac{\beta}{2} + \frac{\beta^2}{4}
       \frac{2}{\beta} \right] =   {n \choose k} \beta,
\end{eqnarray*}
where the last inequality follows from the fact that the function $(0,\infty)\ni x \mapsto x +{{\beta^2}/{4x}}$
attains its minimal value at $x = {{\beta}/2}$. In the case $i+j <0$ it is enough to replace $i, j, \alpha,
\beta$ by $-i, -j, {1/{\alpha}}, {{\beta}/{\alpha^2}}$, respectively, in the previous reasoning. \qed

\vspace{2mm}

\noindent
\begin{lemma}
Let $a \in (0,\infty)$, $b \in (0,1)$ and let $f \colon [0,\infty) \rightarrow \mathbb{R}$ be a solution of
$(1)$ continuous at zero and such that $f(0)=0$. Then for every $k \in \Nn$ and every $t>0$ we have
$$
f(t) = \sum_{n=k}^{\infty} {n\choose k} f\left( a^{k+1} b^{n-k}
t \right).
$$
In particular, for every $k \in \Nn$, and every $t>0$
$$
\lim_{n\rightarrow \infty} {n\choose k} f\left( a^{k} b^{n-k} t
\right) = 0.
$$
\end{lemma}

\vspace{2mm}

\noindent
{\bf Proof.} Let $t>0$. By mathematical induction we have that
for every $m\in \Nn$
$$
f(t) = \sum_{n=0}^{m-1} f(ab^n t) + f(b^m t) \rightarrow
\sum_{n=0}^{\infty} f(ab^n t),
$$
since $b\in(0,1)$ and $f$ is continuous at zero, so $f(b^m t) \rightarrow 0$ for $m \rightarrow \infty$. This
means that we have proved the required equality in the case $k=0$. This equality will be used in the next step
of the proof.

Now let $k\in \Nn$ be fixed and assume that for each $t>0$
$$
f(t) = \sum_{n=k-1}^{\infty} {n\choose {k-1}} f\left( a^{k}
b^{n-k+1} t \right).
$$
Then we have
\begin{eqnarray*}
f(t) & = & \sum_{n=k-1}^{\infty} {n\choose {k-1}} f\left( a^{k}
         b^{n-k+1} t \right) = \sum_{n=k}^{\infty}
         {{n-1}\choose {k-1}} f\left( a^{k} b^{n-k} t\right) \\
  & = & \sum_{n=k}^{\infty} {{n-1}\choose {k-1}}
         \sum_{j=0}^{\infty} f\left( a^{k+1} b^{n+j -k}
         t\right) = \sum_{j=0}^{\infty} \sum_{n=k}^{\infty}
         {{n-1}\choose {k-1}} f\left( a^{k+1} b^{n+j -k} t\right) \\
  & = & \sum_{i=n}^{\infty} \sum_{n=k}^{\infty} {{n-1}\choose
         {k-1}} f\left( a^{k+1} b^{i -k} t\right)  =
         \sum_{i=k}^{\infty} \left[ \sum_{n=k}^{i}
         {{n-1}\choose {k-1}} \right] f\left( a^{k+1} b^{i -k}
         t\right) \\
  & = &  \sum_{i=k}^{\infty} {{i}\choose {k}} f\left( a^{k+1}
         b^{i -k} t\right),
\end{eqnarray*}
which ends the proof. \qed

\vspace{2mm}

\noindent
\begin{lemma}
Let $a,b \in (0,\infty)$ and let $f\in \mathcal{F}$ be a solution of
 $(1)$. Then for every $i,j \in {\mathbb Z}$
$$
f(t)^2 \leqslant f\left( a^{-i} b^{-j} t \right) f\left( a^i b^j t\right), \hspace{5mm} t>0.
$$
\end{lemma}

\vspace{2mm}

\noindent {\bf Proof.} Lemma 2 yields $a,b \in (0,1)$. By Lemma 3 we may assume, without loss of generality,
that $f$ is positive. Let $i,j \in \mathbb{Z}$ and $\alpha>0$, $\beta \in (0, 2\alpha)$ be fixed. Choose $n_0
\in \Nn$ for which the assertion of Lemma 4 holds. Let $t>0$ and $\varepsilon >0$. By Lemma 5 there exists $n
\geqslant 2n_0$ such that
$$
{n\choose k} f\left( a^{k} b^{n-k} t \right) \leqslant \frac{\varepsilon}{2\beta n_0}\hspace{3mm} \hbox{  and  }
\hspace{3mm} {n\choose k} f\left( a^{n-k} b^{k} t \right) \leqslant \frac{\varepsilon}{2 \beta n_0}
$$
for every $k \in \{ 0,1,\dots, n_0 -1\}$. Then
$$
{n\choose k} f\left( a^{k} b^{n-k} t \right) \leqslant \frac{\varepsilon}{2\beta n_0}, \hspace{5mm} k \in
\{0,\dots, n_0 - 1\} \cup \{ n-n_0 +1, \dots, n \}.
$$
Now, by Lemma 1, we obtain
\begin{eqnarray*}
\lefteqn{f\left( a^{-i} b^{-j} t \right) + \alpha^2 f\left(
   a^i b^j t\right)} \\
 && = \sum_{k=0}^{n+i+j} {{n+i+j} \choose{k}} f\left( a^{k-i}
   b^{n+ i -k} t \right) + \alpha^2 \sum_{k=0}^{n-i-j}
   {{n - i - j}\choose {k}} f\left( a^{k+i} b^{n-i -k}
   t\right) \\
 && = \sum_{m=-i}^{n+j} {{n+i+j} \choose{m+i}} f\left( a^{m}
   b^{n- m} t \right) + \alpha^2 \sum_{m=i}^{n-j}
   {{n - i - j}\choose {m-i}} f\left( a^{m} b^{n-m}
   t\right) \\
  && \geqslant \beta \sum_{m=n_0}^{n-n_0} {{n} \choose{m}} f
   \left( a^{m} b^{n- m} t \right) \geqslant  \beta \sum_{m=
   0}^{n} {{n}\choose {m}} f\left( a^{m} b^{n-m}
   t\right) - \varepsilon \\
  && = \beta f(t) - \varepsilon.
 \end{eqnarray*}
Since $\beta$ was an arbitrary number from the interval
$(0,2\alpha)$ we obtain that for every $t>0$
$$
f\left( a^{-i} b^{-j} t \right) + \alpha^2 f\left(
   a^i b^j t\right) \geqslant 2\alpha f(t),
$$
which implies the required inequality. \qed

The proof of Theorem 1 presented below can be essentially shortened by using the main result from \cite{lacz}.
Nevertheless we decided to give an immediate and elementary argument following some ideas of \cite{Jar}.
However, first notice the following fact concerning subgroups of the multiplicative group $(0,\infty)$.

\vspace{2mm}

\noindent {\large\bf Remark 1.} It is well known (cf. \cite{Hardy Wright}, Ch. XXIII, Th. 438) that every
subgroup of $((0,\infty), \cdot)$ is either of the form $\{q^n \colon n \in \mathbb{Z}\}$ with some $q \in
[1,\infty)$, or is a dense subset of $(0,\infty)$. The second case occurs, for instance, when at least two
elements $x,y$ of the group have non-commensurable logarithms: ${{\ln{x}}/{\ln{y}}} \not\in \mathbb{Q}$ or the
group contains elements different from one but arbitrarily close to it. In particular, each subgroup with
non-empty interior is equal to $(0,\infty)$.

\vspace{2mm}

\noindent {\large\bf Proof of Theorem 1.} By Lemma 3  we may assume, without loss of generality, that $f(t) > 0$
for every $t>0$. It follows from Lemma 2 that $a,b \in (0,1)$. For every $t > 0$ we define the orbit $C(t)$ of
$t$ by
$$
C(t) = \left\{ a^i b^j t \colon\, i,j \in \mathbb{Z} \right\}.
$$
Since $C(1)$ is a subgroup of the group $(0,\infty)$, Remark 1 yields that either ${{\ln{a}}/{\ln{b}}} \in
\mathbb{Q}$ and $C(t)= \{ c^j t \colon j \in \mathbb{Z}\}$ for some $c
>1$, or ${{\ln{a}}/{\ln{b}}} \not\in \mathbb{Q}$ and the set $C(t)$ is dense in $(0,\infty)$.

Take $t_0>0$ and let $G = \overline{\{ \ln{x}\colon x\in C(1)\}}$ and $A = \overline{\{ \ln{x}\colon x\in
C(t_0)\}}$. The function $F\colon A \mapsto \Rr$, defined by $F(x) = \ln {f(e^x)}$, is continuous and, by Lemma
6, convex:
$$
2F(x) \leqslant F(x-h) + F(x+h), \hspace{5mm} x\in A, \, h \in G.
$$
Therefore, if $G=\Rr$ then, according to \cite{Ku} VII.3, Th.
2, the function
$$
A \times \left(G \setminus \{0\} \right) \, \ni (x,h) \mapsto
\frac{F(x+h) - F(x)}{h}
$$
is increasing with respect to each variable. Clearly the same holds
true if the group $G$ is discrete. In particular, we obtain that
the following limits exist and do not depend on $h$:
$$
p_{-} = \lim_{x\rightarrow - \infty} \frac{F(x+h) - F(x)}{h}, \,\,\, \hbox{ and } \,\,\, p_{+} =
\lim_{x\rightarrow \infty} \frac{F(x+h) - F(x)}{h}.
$$
Obviously $-\infty \leqslant p_{-} \leqslant p_{+} \leqslant \infty$. If $p_{+} = \infty$, then, taking $x =
\ln{t}$, $h=\ln{a}$ and $h=\ln{b}$, we would have
$$
\lim_{t\rightarrow \infty} \frac{f(at)}{f(t)} = 0, \,\, \hbox{
and } \,\, \lim_{t\rightarrow \infty} \frac{f(bt)}{f(t)} = 0,
$$
which is impossible in view of equation $(1)$. Thus $p_{+} < \infty$. Similarly $p_{-} > - \infty$. Now we have
that
$$
\lim_{t\rightarrow 0} \frac{f(at)}{f(t)} = a^{p_{-}}, \,\, \hspace{5mm} \,\, \lim_{t\rightarrow \infty}
\frac{f(at)}{f(t)} = a^{p_{+}},
$$
and
$$
\lim_{t\rightarrow 0} \frac{f(bt)}{f(t)} = b^{p_{-}}, \,\, \hspace{5mm} \,\, \lim_{t\rightarrow \infty}
\frac{f(bt)}{f(t)} = b^{p_{+}}.
$$
The equation $(1)$ implies that $a^{p_{-}} + b^{p_{-}} = 1$ and $a^{p_{+}} + b^{p_{+}} = 1$. As $a,b \in (0,1)$
there is exactly one real root $p$ of the equation $a^{p} + b^{p} = 1$; clearly $p > 0$. Consequently, we get
$p_{-} = p_{+} = p$.

Notice now that the monotonicity of $A \ni x \mapsto h^{-1}(F(x+h) - F(x))$ implies the monotonicity of the
continuous functions $C(t_0) \ni t \mapsto {{f(at)}/{f(t)}}$ and $C(t_0) \ni t \mapsto {{f(bt)}/{f(t)}}$, and
thus
$$
 \hspace{5mm} f(at) = a^{p} f(t), \,\, \hbox{ and }\,\,
f(bt) = b^{p}f(t), \hspace{5mm} t \in \overline{C(t_0)} \setminus \{0\}.
$$
If $t\in C(t_0)$ then $t = a^i b^j t_0$ for some $i,j \in
\mathbb{Z}$, so
$$
f(t) = f(a^i b^j t_0) = a^{p i} b^{p j} f(t_0) = \frac{f(t_0)}{t_0^{p}} t^{p}.
$$
This implies that the positive and continuous function $H$ on $(0,\infty)$, defined by
$$
H(t) = t^{-p} f(t),
$$
is constant on every orbit, and thus, in particular, $H(t) =
H(at) = H(bt)$ for each $t \in (0,\infty)$. To complete the
proof it is enough to observe that if ${{\ln{a}}/{\ln{b}}}
\not\in \mathbb{Q}$ then $\overline{C(t_0)} = [0,\infty)$ for
each $t_0 >0$, and thus $H$ is constant on $[0,\infty)$. \qed

\vspace{2mm}

Observe that the property of $H$ in Theorem 1 is similar to that of double periodicity in complex function
theory, which goes back to Jacobi and Weierstrass (see eg. \cite{Copson}, Ch. XIII, XIV).

\section{Stable distributions in the sense of the weak generalized
convolution}

In \cite{Urbanik64} K. Urbanik considered stable distributions with respect to a generalized convolution
$\diamond$. Theorem 4 in \cite{Urbanik64} states that if $h$ is a continuous homomorphism from $\mathcal{P}_{+}$
to $\mathbb{R}$ and for some measure $\lambda \in \mathcal{P}_{+}$ there exists a sequence
$(c_n)_{n\in\mathbb{N}}$ of positive numbers such that
$$
T_{c_n} \lambda^{\diamond n} \rightarrow \lambda_0 \neq \delta_0,
$$
then there exist positive numbers $c$ and $p$ such that
$$
h(T_x \lambda_0) = \exp\left\{ - c x^p \right\}, \hspace{5mm}  x\geqslant 0.
$$
The existence of at least one such measure $\lambda$ follows from the property (v) of the definition of the
generalized convolution $\diamond$; thus we know that for $\lambda = \delta_1$ such a sequence $(c_n)_{n\in
\mathbb{N}}$ exists. Because of the form of $h(T_x \lambda_0)$ it was natural to call the measure $\lambda_0$
stable with respect to the convolution $\diamond$.

For the weak generalized convolution $\oplus_{\mu}$ the condition (v) need not hold, so we will define
$\oplus_{\mu}$-stable distributions using the classical linearity conditions. In this paper we restrict our
attention to the symmetric weakly stable distributions.

\begin{definition}
Let $\mu$ be a non-trivial symmetric weakly stable measure on a separable Banach space $\mathbb{E}$. A measure
$\lambda \in {\cal P}_{+}$ is stable with respect to the weak generalized convolution $\oplus = \oplus_{\mu}$ if
$$
\forall \,\, r,s\geqslant 0 \,\, \exists \, c(r,s)\geqslant 0, d(r,s) \in \mathbb{R} \hspace{5mm}  (T_r \lambda)
\oplus (T_s \lambda) = (T_{c(r,s)} \lambda) \oplus \delta_{d(r,s)}.
$$
If for every $r,s>0$ we have $d(r,s) = 0$ then we say that
$\lambda$ is strictly stable with respect to $\oplus$.
\end{definition}

\noindent {\large\bf Remark 2.} The condition in Definition 2 written in the language of random variables states
that there are functions $c \colon [0,\infty)^2 \rightarrow [0,\infty)$ and $d \colon [0,\infty)^2 \rightarrow
\mathbb{R}$ such that
$$
 \forall \,\, r,s\geqslant 0 \hspace{5mm}  r \mathbf{X} \theta + s \mathbf{X}' \theta'
\stackrel{d}{=} c(r,s) \mathbf{X} \theta + d(r,s) \mathbf{X}',
$$
where $\mathcal{L}(\mathbf{X}) = \mathcal{L}(\mathbf{X}') = \mu$, $\mathcal{L}(\theta) = \mathcal{L}(\theta') =
\lambda$, $\mathbf{X}, \mathbf{X}', \theta, \theta'$ are independent. Putting $\|(r,s)\|_2 = \sqrt{r^2 + s^2}$
for every $r,s\geqslant 0$, $(r,s) \neq (0,0)$, we have
\begin{eqnarray*}
\lefteqn{ \frac{r}{\|(r,s)\|_2} \mathbf{X} \theta + \frac{s}{\|(r,s)\|_2} \mathbf{X }' \theta'} \\
    & \stackrel{d}{=} & c\left(\frac{r}{\|(r,s)\|_2},\frac{s}{\|(r,s)\|_2}\right) \mathbf{X} \theta +
d\left(\frac{r}{\|(r,s)\|_2}, \frac{s}{\|(r,s)\|_2}\right) \mathbf{X}',
\end{eqnarray*}
whence
$$
\forall \,\, r,s \geqslant 0 \hspace{5mm} r \mathbf{X} \theta + s \mathbf{X}' \theta' \stackrel{d}{=} c_1(r,s)
\mathbf{X} \theta + d_1(r,s) \mathbf{X}',
$$
where $c_1$ and $d_1$ are given by $c_1(0,0) = d_1(0,0)= 0$
$$
c_1(r,s) = \|(r,s)\|_2\ c\left( \frac{r}{\|(r,s)\|_2}, \frac{s}{\|(r,s)\|_2}\right),
$$
and
$$
d_1(r,s) = \|(r,s)\|_2\ d\left( \frac{r}{\|(r,s)\|_2}, \frac{s}{\|(r,s)\|_2}\right).
$$
In the other words, in Definition 2 we can always assume that the functions $c$ and $d$ are homogeneous:
$$
c(ru, su) = uc(r,s), \hspace{4mm} \hbox{and} \hspace{4mm} d(ru, su) = ud(r,s),
$$
for every $r,s,u \geqslant 0$. Notice also that, since $\mathbf{X}$ has a symmetric distribution, $a\mathbf{X}
\stackrel{d}{=} -a\mathbf{X}$ for every $a \in \mathbb{R}$. This implies that without loss of generality we can
assume that $d$ is nonnegative, taking if necessary $| d(r,s)|$ instead of $d(r,s)$ for every $r,s \geqslant 0$.
In what follows we will always do so. On the other hand, as we will see in Remark 4, the functions $c$ and $d$,
at least in the case of symmetric stable distribution $\mu = \gamma_{p}$, need be neither
homogeneous, nor unique, nor continuous. \\

Notice that for every $x\in \mathbb{R}$ the measure $\delta_x$ is stable with respect to the weak generalized
convolution $\oplus = \oplus_{\mu}$ for each weakly stable measure $\mu$. Indeed, for every $r,s >0$ we have
$(T_r \delta_x) \oplus (T_s \delta_x) = (T_r \delta_x) \oplus \delta_{sx}$, so the condition holds with
$c(r,s) = r$ and $d(r,s) = sx$. Consequently, we will say that $\delta_x$ is the trivial example of
$\oplus_{\mu}$-stable distribution, just as in the classical case.

In order to characterize stable distributions with respect to the weak generalized convolution $\oplus_{\mu}$ we
choose first $\xi \in \mathbb{E}^{\ast}$ such that the random variable $<\xi, \mathbf{X}>$,
$\mathcal{L}(\mathbf{X}) = \mu$, is symmetric and non-trivial. Since the distribution $\mu_{\xi} =
\mathcal{L}(<\xi, \mathbf{X}>)$ is weakly stable for each weakly stable measure $\mu$, we have that if
$\mu_{\xi} \circ \lambda_1 = \mu_{\xi} \circ \lambda_2$ for some $\lambda_1, \lambda_2 \in \mathcal{P}_{+}$ then
$\lambda_1 = \lambda_2$.

Following the Urbanik construction for every $\xi \in \mathbb{E}^{\ast}$ we define a homomorphism $h_{\xi}
\colon \mathcal{P}_{+} \rightarrow \mathbb{R}$ by the formula
$$
h_{\xi}(\nu) = \int_{[0,\infty)} \widehat{\mu}(t\xi) \nu(dt), \hspace{6mm} \nu \in \mathcal{P}_{+}.
$$
Let $\psi(t) = \widehat{\mu}(t\xi)$, $t \in \mathbb{R}$. We see that $\psi$ is the characteristic function of
the measure $\mu_{\xi}$ and the function $\mathbb{R} \ni t \mapsto h_{\xi}(T_t \nu)$ is the characteristic
function of the measure $\mu_{\xi} \circ \nu$. From the previous considerations it follows that the functionals
$h_{\xi}(T_t \cdot)$, $t \in \mathbb{R}$, separate points in $\mathcal{P}_{+}$ in the sense that if $h_{\xi}(T_t
\lambda_1) = h_{\xi}(T_t \lambda_2)$ for some $\lambda_1, \lambda_2 \in \mathcal{P}_{+}$ and for each $\xi \in
\mathbb{E}^{\ast}$ and $t \in \mathbb{R}$, then $\lambda_1 = \lambda_2$. \\

Now for every $t \in \mathbb{R}$ let $\varphi(t) = h_{\xi}(T_t \lambda)$, where $\lambda \in \mathcal{P}$ is
stable with respect to the generalized convolution $\oplus_{\mu}$. By Definition 2  there exist nonnegative
functions $c$ and $d$ on $[0,\infty)^2$ such that the functional equation
$$
\varphi(rt) \varphi(st) = \varphi(c(r,s)t) \psi(d(r,s)t) \eqno{(2)}
$$
is satisfied. According to Remark 2 we can also assume, if necessary, that $c$ and $d$ are homogenous. We need
to solve equation (2) in the set $\Phi_2$ of pairs of characteristic functions, given by
$$
\Phi_2 = \left\{ (\varphi, \psi) \in \Phi \times \Phi \colon \exists \,\, \lambda \in \mathcal{P}_{+} \,\,\,
\forall \,\, t \in \mathbb{R}\hspace{3mm}\varphi(t) = \int_{\mathbb{R}} \psi(ts) \lambda(ds) \right\}.
$$
In what follows we consider four possible, complementary but not disjoint, cases concerning the function $d$:
\begin{namelist}{ll}
\item[3.1] measures strictly stable with respect to the weak generalized convolution, when $d (r,s) = 0$ for
           every $r,s>0$;
\item[3.2] $p$-self-decomposable measures with respect to the weak generalized convolution, when
           such that $d(r,0)> 0$ or $d(0,r) >0$ for some $r>0$;
\item[3.3] semi-stable measures with respect to the weak generalized convolution, when $d(r,s)=0$ for some
           $r,s >0$;
\item[3.4] $(c,d)$-pseudostable measures with respect to the weak generalized convolution, when $d(r,s)
           > 0$ for each $r,s >0$.
\end{namelist}

\subsection{Measures strictly stable with respect to the weak
generalized convolution}

Definition 2 states that if $d(r,s)=0$, for every $r,s \geqslant 0$, then the measure $\lambda$ is strictly
stable with respect to $\oplus_{\mu}$. In this case we do not need to assume that the measure $\mu$ is
symmetric, so (2) leads to the equation
$$
\forall \,\, r,s\geqslant 0 \,\, \exists \,\, c(r,s)\geqslant 0 \,\, \forall t \in \Rr \hspace{5mm} \varphi(rt)
\varphi(st) = \varphi(c(r,s)t). \eqno{(3)}
$$
Equation (3) is the classical one defining strictly stable distributions and neither symmetry of the measure
$\mu$, nor homogeneity of the function $c$ is required to obtain that $\varphi$ is the  characteristic function
of a strictly stable distribution. This means that there exist $p\in (0,2]$, $\sigma>0$ and $\beta \in [-1, 1]$
such that
$$
\varphi(t) = \left\{ \begin{array}{ll}
    \exp\left\{ - \sigma^p |t|^p \left(1- \beta \hbox{sgn}(t) \tan{ \frac{\pi p}{2}} \right) \right\} & \hbox{if
    }\ p \neq 1, \\
     \exp\left\{ - \sigma |t| + i\mu t \right\} & \hbox{if } \ p = 1.
\end{array} \right.
$$

 Since $\varphi$ is the characteristic function of the random variable $<\xi, {\bf X} \theta$>, where
${\cal L}({\bf X}) = \mu$, ${\cal L}(\theta) = \lambda$ , ${\bf X}$ and $\theta$ independent, then we obtain
that every one-dimensional projection of the random vector ${\bf X} \theta$ is strictly stable. This implies
(for details see e.g. Th. 2.1.5 in \cite{ST}) that the index of stability does not depend on $\xi$ and the
random vector ${\bf X} \theta$ is strictly stable. In this way we proved the following:

\begin{theorem}
Let $\mu$ be a non-trivial weakly stable distribution on $\mathbb{E}$. If $\lambda \in \mathcal{P}$ is strictly
stable with respect to the weak generalized convolution $\oplus_{\mu}$ then $\mu \circ \lambda$ is strictly
stable, i.e. $\mu$ and $\lambda$ are factors of a strictly stable distribution.
\end{theorem}

\vspace{2mm}

\noindent {\large\bf Example 1.} Let ${\bf U}^n$ be the random vector with the uniform distribution $\mu =
\omega_n$ on the unit sphere $S_{n-1}\subset \mathbb{R}^n$ (as well we can consider here any projection ${\bf
U}^{n,k}$ of ${\bf U}^n$ into $\mathbb{R}^k$, $k< n$). The $\omega_n$-weakly Gaussian random variable $\Gamma_n$
is defined by the following equation:
$$
{\bf U}^n \cdot \Gamma_n \stackrel{d}{=} (X_1,\dots,X_n)= {\bf {\bf X}},
$$
where ${\bf U}^n$ and $\Gamma_n$ are independent, ${\bf X}$ is an $n$-dimensional Gaussian random vector with
independent identically distributed coordinates. It is known, and it was already known to Schoenberg in 1938
(see \cite{Schoenberg1}), that for every spherically invariant random vector $\mathbf{Y}$ in $\mathbb{R}^n$ we
have $\mathbf{Y} \stackrel{d}{=} {\bf U}^n \cdot \| \mathbf{Y} \|_2$, where ${\bf U}^n$ and $\|\mathbf{Y}\|_2$
are independent. This implies that $\Gamma_n$ is the distribution of $\|{\bf X}\|_2$. Simple calculations show
that $\Gamma_n$ has the density given by
$$
f_{2,n} (r) = \frac{2}{2^{n/2} \Gamma(\frac{n}{2})} r^{n-1} e^{-{{r^2}/2}}.
$$
For $n=2$ this is the Rayleigh distribution with parameter $\lambda = 2$, thus the Rayleigh distribution is
$\omega_2$-weakly Gaussian. For $n=3$ this is the Maxwell distribution with parameter $\lambda = 2$, thus the
Maxwell distribution is $\omega_3$-weakly Gaussian. Recall that the generalized Gamma distribution with
parameters $\lambda, p, a >0$ (notation $\Gamma(\lambda, p,a)$) has density function given by
$$
f(x) = \frac{a}{\Gamma({p/a}) \lambda^{{p/a}}} x^{p-1} \exp \left\{ - \frac{x^a}{\lambda} \right\}, \hspace{7mm}
x>0.
$$
Thus we have that the generalized Gamma distribution $\Gamma(\lambda, n, 2)$ is
$\omega_n$-weakly Gaussian.

Now let $\theta_{\alpha}^n$ be an $\omega_n$-weakly strictly $\alpha$-stable random variable. Then ${\bf U}^n
\theta_{\alpha}^n$ is rotationally invariant $\alpha$-stable random vector for the vector ${\bf U}^n$
independent of $\theta_{\alpha}^n$. On the other hand every rotationally invariant $\alpha$-stable random vector
has the same distribution as ${\bf Y} \sqrt{\theta_{{\alpha}/2}}$, where ${\bf Y}$ is a rotationally invariant
Gaussian random vector independent of the nonnegative variable $\theta_{{\alpha}/2}$ with the Laplace transform
$e^{-t^{{\alpha}/2}}$. Finally we have
$$
{\bf U}^n \cdot \theta_{\alpha}^n \stackrel{d}{=} {\bf U}^n \Gamma_n \sqrt{\theta_{{\alpha}/2}},
$$
for ${\bf U}^n$, $\Gamma_n$ and $\theta_{{\alpha}/2}$ independent. This implies that the density of a
$\omega_n$-weakly strictly $\alpha$-stable random variable $\theta_{\alpha}^n$ is given by
$$
f_{\alpha,n} (r) = \int_0^{\infty} f_{2,n}\left(\frac{r}{\sqrt{s}}\right)
\frac{1}{\sqrt{s}} f_{{\alpha}/2}(s) ds.
$$
In particular, if we take $\alpha = 1$ then
$$
f_{{1/2}} (s) = \frac{1}{\sqrt{2\pi}} x^{-{3/2}} e^{- {1/{(2x)}}}, \hspace{5mm} x>0.
$$
Simple calculations and the duplication formula
$$
\sqrt{\pi} \Gamma(2s) = 2^{2s -1} \Gamma(s) \Gamma\left(s + \frac{1}{2}\right), \hspace{5mm} s>0,
$$
show that
$$
f_{1,n} (r) = \frac{2^{2-n} \Gamma(n)}{\Gamma({n/2}) \Gamma({n/2})} \frac{r^{n-1}}{(r^2 + 1 )^{{(n+1)}/2}}, \hspace{5mm} r>0,
$$
is the density function of the $\omega_n$-weakly strictly Cauchy distribution.

\vspace{2mm}

\noindent {\large\bf Remark 3.} Finally we get the following interesting statement:

{\em If $X$ is a nonnegative random variable with the density function $f_{\alpha,n}$, or if $X$ is a random
variable such that $\mathcal{L}(|X|)$ has the density $f_{\alpha, n}$ and $X'$ has also this property, then
the random variable
$$
Y = \left\| X \mathbf{U}^n + X' \mathbf{U}'^n \right\|_2
$$
has the density $2^{-1/{\alpha}} f_{\alpha,n}\left(\mathbf{x} \cdot 2^{-1/{\alpha}}\right)$, i.e. $Y
\stackrel{d}{=} 2^{1/{\alpha}} |X|$.}

\vsp

\noindent {\large\bf Remark 4.} Spherically invariant (or spherically generated, or rotationally invariant)
measures, mentioned in Example 1, are extensively studied and applied in stochastic modeling. More information
about such measures can be found in \cite{FKotzNg}. It is worth mentioning here that in 1963 J.F.C. Kingman (see
 \cite{Kingman}) constructed an independent increments, two-dimensional (in the simplest case) stochastic process
$\{ X_t\colon t \geqslant 0\}$, where increments are spherically invariant. This process is associated to $\{Y_t
= \| X_t \|_2 \colon t \geqslant 0\}$, describing the distance of the particle to the origin, has increments
independent in the sense of the $\omega_2$-weak generalized convolution, but obviously, these increments are not
independent in the usual sense. This paper was an important part of the original motivation for Urbanik's
generalized convolution.

\vsp

The next theorem gives the full characterization of $\gamma_{\alpha}$-weakly strictly stable distribution
$\lambda$ for $\gamma_{\alpha}$ strictly $\alpha$-stable. Recall that among  stable distributions only strictly
stable distributions are weakly stable and, except the symmetric case, only $\mathbb{R}_{+}$-weak stability can
be considered, i.e. constants $a,b,c$ in definition of weakly stable distribution shall be positive (for details
see \cite{MOU}). Non-symmetric stable distributions are not weakly stable on $\Rr$.

\begin{theorem}
Let $\gamma_{\alpha}$ be a strictly $\alpha$-stable distribution on $\Rr^n$. Then the following conditions are
equivalent:
\begin{itemize}
\item[{\rm 1.}] $\lambda$ is $\gamma_{\alpha}$-weakly strictly
          stable;
\item[{\rm 2.}] $\lambda = \mathcal{L}(a \theta_p^{1/{\alpha}})$ for some $a>0$, $p\in (0,1]$ and $\theta_p$ is
          a positive random variable with  Laplace transform $\Ee e^{-t \theta_p} = e^{-t^p}$ if $p\in
          (0,1)$, and $\theta_1 \equiv 1$;
\item[{\rm 3.}] $\gamma_{\alpha} \circ \lambda = T_a \gamma_{\alpha
          p}$ for some $\alpha p$-strictly stable distribution
          $\gamma_{\alpha p}$ and some numbers $a>0$ and $p\in (0,1)$.
\end{itemize}
\end{theorem}

\vspace{2mm}

\noindent {\large\bf Proof.} By Theorem 2 we have that condition 1 implies that the scale mixture
$\gamma_{\alpha} \circ \lambda$ is strictly stable; thus $1$ implies $3$. Evidently $2$ yields $1$. The
implication $3 \Rightarrow 2$ follows from \cite{MBO2}, where M. Borowiecka-Olszewska proved that the scale
mixture of strictly $\alpha$-stable distribution is stable if and only if the mixing distribution $\lambda$
equals $\mathcal{L}(a \theta_p^{1/{\alpha}})$ for some $a>0$ and $p \in (0,1]$.  \qed

\subsection{c-selfdecomposable measures with respect to the\\ weak generalized convolution}

The case which we are considering here seems to be rather unrealistic but it leads to a very interesting class
of distributions and, because of this, is worth including.

Assume, for instance, that for some $r>0$ we have $ d(r,0) \neq 0$, and so, by the homogeneity
condition $d(1,0) \neq 0$. If now $c(1,0) = 0$, then we would have $\varphi(t) = \psi(d(1,0) t)$, $t \in
\mathbb{R}$, and the corresponding measure $\lambda$ would be concentrated at a point. This case is considered
as trivial. Thus assume also that $c(1,0) \neq 0$. Now, putting $\alpha = c(0,1)$ and $\beta = d(0,1)$, we see
that $\alpha, \beta >0$ and $(\varphi, \psi )$ is a solution of the functional equation
$$
\varphi(t) = \varphi (\alpha t) \psi(\beta t).
$$
This recalls the definition of $b$-semi-selfdecomposable distributions, which in the book of Sato \cite{Sato}
was given in the following way:

{\em A non-trivial probability measure $\nu$ is semi-selfdecomposable if there exist a number $b>1$ and an
infinitely divisible  probability measure $\rho_b$ such that
$$
\widehat{\nu}(t) = \widehat{\nu}(b^{-1} t) \widehat{\rho_{b}}(t), \hspace{5mm} t \in \mathbb{R}.
$$}

\vspace{2mm}

Notice that the condition $b>1$ can be omitted if we are talking about nontrivial distributions. If $b=1$, then
the condition holds for any $\nu$ with the trivial measure $\rho_1 = \delta_0$. For $b<1$, notice  first that
\begin{eqnarray*}
\widehat{\nu}(t) & = & \widehat{\nu}(b^{-1} t)
           \widehat{\rho_{b}}(t) = \widehat{\nu}(b^{-2} t)
         \widehat{\rho_{b}}(b^{-1}t) \widehat{\rho_{b}}(t) \\
   & = & \dots = \widehat{\nu}(b^{-n} t) \prod_{j=0}^{n-1}
         \widehat{\rho_{b}}(b^{-j}t).
\end{eqnarray*}
Substituting $t \rightarrow b^n t$, we obtain
$$
\widehat{\nu}(b^n t) = \widehat{\nu}(t) \prod_{j=0}^{n-1}
         \widehat{\rho_{b}}(b^{n-j}t).
$$
The left hand side of this formula tends to 1 when $n\rightarrow \infty$; thus there exists also the
corresponding limit of the right hand side. Since absolute values of all the functions here are bounded from
above by 1, this equality holds only if for every $t\in \mathbb{R}$
$$
\left| \widehat{\nu}(t) \right| = \lim_{n\rightarrow \infty} \left| \prod_{j=0}^{n-1}
\widehat{\rho_{b}}(b^{n-j}t) \right| = 1,
$$
which is impossible as $\nu$ is non-trivial.

In our case $\widehat{\nu} = \varphi$, $b^{-1} = \alpha$ and $\widehat{\rho_b} = \psi ( \beta \cdot)$, but we
cannot assume that $\psi$ is an infinitely divisible characteristic function. All we know is that
$\psi$ is the characteristic function of some nontrivial weakly stable distribution. Notice that
$$
\varphi( t) = \varphi(\alpha^n t) \prod_{j=0}^{n-1} \psi(\beta \alpha^{j}t) = \prod_{j=0}^{\infty} \psi(\beta
\alpha^{j}t),
$$
since we have already shown that $\alpha<1$. Now we see that $\varphi$ is the characteristic function of the
random variable
$$
Y = \beta \sum_{j=0}^{\infty} \alpha^j X_j,
$$
where $X_j$, $j=0,1,\dots$, are independent copies of the variable $X$ with distribution $\mu$. Since the class
$\{ \mu \circ \lambda \colon \lambda \in \mathcal{P}\}$ of mixtures of the weakly stable distribution $\mu$
forms a set which is weakly closed and closed under convolution and rescaling, we see that $Y \stackrel{d}{=} X
\theta$ for some $\theta$ independent of $X$, $\mathcal{L}(\theta) = \lambda$, under the assumption that the
series defining $Y$ converges in distribution. Now let us accept

\begin{definition}
Let $X$ be a non-trivial weakly stable random vector  with the distribution  $\mu$ and let $\alpha \in (0,1)$. A
random variable $\theta$ (or its distribution) is $\alpha^{-1}$-selfdecomposable in the sense of the generalized
weak convolution $\oplus_{\mu}$ if there exists a random variable $Q$ such that
$$
\theta = \left( \alpha \theta \right) \oplus_{\mu} Q.
$$
\end{definition}

The previous considerations show that the following proposition holds:

\begin{prop}
Let $\mathbf{X}$ be a non-trivial weakly stable random vector  with the distribution  $\mu$, $\alpha \in(0,1)$
and let $\theta$ be $\alpha^{-1}$-selfdecomposable in the sense of the weak generalized convolution
$\oplus_{\mu}$. Assume that $\mathbf{X}$ and $\theta$ are independent. Then
\begin{itemize}
\item if $\theta = \left( \alpha \theta \right) \oplus_{\mu} \beta$, then the random variable $\theta$ is
uniquely determined by the condition
$$
\mathbf{X} \theta \stackrel{d}{=} \beta \sum_{j=0}^{\infty} \alpha^j \mathbf{X}_j,
$$
where $X_j$'s are independent copies of $X$. %
\item if $\theta = \left( \alpha \theta \right) \oplus_{\mu} Q$ for some non-trivial random variable $Q$, then
$\theta$ is uniquely determined by the condition
$$
\mathbf{X} \theta \stackrel{d}{=} \sum_{j=0}^{\infty} \alpha^j \mathbf{X}_j\  Q_j,
$$
where $\mathbf{X}_j$'s are independent copies of $\mathbf{X}$ and $Q_j$'s are independent copies of $Q$ such
that $(\mathbf{X}_j)_{j\in \mathbf{N}}$ and $(Q_j)_{j\in \mathbf{N}}$ are independent.
\end{itemize}
\end{prop}

\vspace{1mm}

 In the case when $X$ has strictly $p$-stable distribution, we have that the series defining $Y$ converges at
 least in distribution, and
$$
Y = \beta \sum_{j=0}^{\infty} \alpha^j X_j \stackrel{d}{=} \Biggl( \sum_{j=0}^{\infty} \alpha^{jp} \Biggr)^{1/p}
X = \left( 1- \alpha^p \right)^{-{1/p}} X.
$$
Thus the random variable $\theta$ exists and $\Pp \{ \theta = \left( 1- \alpha^p \right)^{-{1/p}} \} = 1$, which
we consider as a trivial solution.

Let $X = U_{n,1}$ with the distribution $\omega_{n,1}$ be the one-dimensional margin of the vector
$\mathbf{U}^{n}$ with the uniform distribution on the unit sphere in $\mathbb{R}^n$. Then we have $\Ee U_{n,1} =
0$ and ${\rm Var} U_{n,1} \leqslant 1 < \infty$. Thus the series defining $Y$ converges in $L_2$, so also almost
everywhere. Thus the random variable $Y$ is well defined, and there exists $\theta_n$ such that
$$
X \theta_n \stackrel{d}{=} \beta \sum_{j=0}^{\infty} \alpha^j X_j.
$$
In the case $n=3$ the variable $U_{3,1}$ has the uniform distribution on the interval $[-1,1]$, thus
$$
\varphi(t) = \prod_{j=0}^{\infty} \frac{ \sin{(\beta \alpha^j t)}}{\beta \alpha^j t}.
$$
Let $G$ be the cumulative distribution function for $Y = U_{3,1} \theta_3$. Notice that for every $u \in [0,
(1-\alpha)^{-1}]$
$$
G(u) = \Pp \left\{ X_1 + \alpha \left( X_2 + \alpha X_3 + \dots \right) < {u/{\beta}} \right\} = \Pp \left\{ X_1
+ \alpha Y' < {u/{\beta}} \right\},
$$
where $Y'$ is a copy of $Y$ independent of $X_1$. One can also show that for every $u \in [-
\beta(1-\alpha)^{-1}, \beta (1-\alpha)^{-1}]$ we have
$$
\beta G'(u) =  G \left( \frac{u+\beta}{\alpha} \right) - G \left( \frac{u-\beta}{\alpha} \right).
$$
Moreover, $G$ as the distribution function of a symmetric random variable has the property $G(-u) = 1 - G(u)$
for every $u \in \mathbb{R}$. $G$ is also strictly increasing on $[-\beta(1-\alpha)^{-1}, \beta
(1-\alpha)^{-1}]$ and $G (-\beta(1-\alpha)^{-1})= 0$. Similar conditions can be obtained considering $U_{n,1}$
for an arbitrary $n \in \mathbb{N}$.

\subsection{Semi-stable measures with respect to the weak generalized convolution}

Let $Z = \left\{ (r,s) \in (0,\infty)^2 \colon d(r,s) =0\right\}$. In this subsection we assume that $Z$ is
nonempty, postponing the case $Z = \emptyset$ to the next section. We assume also that $Z \neq (0,\infty)^2$,
since the opposite case was considered in section 3.1.

\vspace{2mm}

\begin{theorem}
Let $c, d \colon [0,\infty)^2 \rightarrow [0,\infty)$. Assume that $Z \neq \emptyset$ and $Z \neq (0,\infty)^2$.
Let $\varphi,\psi \colon \mathbb{R} \rightarrow \mathbb{R}$ be nontrivial characteristic functions satisfying
equation $(2)$. Then there exist $p \in (0,2]$ and even continuous functions $H,K \colon \mathbb{R} \setminus \{
0 \} \rightarrow (0,\infty)$ such that
$$
\varphi(t) = e^{-|t|^{p} H(t)}\hspace{4mm} \hbox{and} \hspace{4mm} \psi(t) = e^{-|t|^{p} K(t)}
$$
for every $t \in \mathbb{R} \setminus \{0\}$ and
$$
H(rt) = H(st) = H(c(r,s)t) \hspace{4mm} \hbox{and} \hspace{4mm} K(rt) = K(st) = K(c(r,s)t)
$$
for every $t \in \mathbb{R} \setminus \{0\}$ and $(r,s) \in Z$. Moreover,
$$
r^{p} H(rt) + s^{p} H(st) = c(r,s)^{p} H(c(r,s)t) + d(r,s)^{p} K(d(r,s)t)
$$
for every $t \in \mathbb{R} \setminus \{0\}$ and $r,s>0$ such that $c(r,s)>0$ and $(r,s) \not\in Z$.
\end{theorem}

\vspace{2mm}

\noindent {\bf Proof.} Take any $(r,s)\in Z$. Then for every $t\in \mathbb{R}$ we have that $\varphi(rt)
\varphi(st) = \varphi\left( c(r,s)t\right)$. The equality $c(r,s) = 0$ would imply that $\varphi(rt) \varphi(st)
\equiv 1$ and next $\varphi \equiv 1$. Thus $c(r,s)>0$, whence $\varphi$ satisfies the functional equation
$$
\varphi(t) = \varphi(at) \varphi( bt), \eqno{(4)}
$$
where $a:= {r/{c(r,s)}}$ and $b := {s/{c(r,s)}}$ are positive numbers. This implies, in particular, that
$\varphi$ is the characteristic function of an infinitely divisible distribution and, consequently, it does not
attain value zero. This property can easily be proved directly. First we verify  that $a,b<1$. Without loss of
generality we can assume that $a\geqslant b$. If $a=1$, then by ($4$), $\varphi(bt) = 1$ for every $t>0$, so
$\varphi$ would be identically one. Suppose that $a>1$ and take any $u>0$ such that $\varphi(t) >0$ for every $t
\in (0,u)$. If $t \in (0,u)$, then $a^{-1}t, a^{-1} bt \in (0,u)$ and thus by ($4$),
$$
1\geqslant \varphi(t) = \varphi(a a^{-1} t) = \frac{\varphi(a^{-1} t)}{\varphi(b a^{-1}t)}\geqslant
\varphi(a^{-1}t),
$$
whence, by induction,
$$
1 \geqslant \varphi(t)\geqslant \varphi(a^{-n} t), \hspace{5mm} n \in \mathbb{N}.
$$
Now, using the continuity of $\varphi$, we deduce that $\varphi(t)=1$ for every $t\in (0,u)$, and, consequently,
for every $t>0$. Therefore $a<1$ and $b<1$. Suppose now that there exists $t_0 >0$ such that $\varphi(t_0) =0$.
Then ($4$) gives $\varphi(t_n) = 0$ for each $n\in \mathbb{N}$, where
$$
t_{n+1} = at_n \hspace{4mm} \hbox{or} \hspace{4mm} t_{n+1} = bt_n, \hspace{5mm} n \in \mathbb{N}.
$$
Since $a,b<1$ the sequence $(t_n)_{n \in \mathbb{N}}$ tends to zero and thus, by the continuity of $\varphi$, we
get $\varphi(0) =0$, which is impossible.

Now we can define $f = -\ln \varphi$. The function $f$ is non-negative, continuous and even, $f(0) = 0$.
Moreover, $f$ as a logarithm of characteristic function does not attain value zero in a vicinity of zero.
Clearly, $f$ satisfies equation $(1)$. By Theorem 1 there exist  $p \in (0,\infty)$ and an even continuous
function $H \colon \mathbb{R}\setminus \{0\} \rightarrow (0,\infty)$  such that for every $t\in
\mathbb{R}\setminus \{0\}$
$$
f(t) = |t|^{p} H(t) \hspace{4mm} \hbox{and} \hspace{4mm} H(rt) = H(st) = H(c(r,s)t).
$$
Consequently,
$$
\varphi(t) = e^{-|t|^{p} H(t)}, \hspace{5mm} t \in \mathbb{R}\setminus \{0\}.
$$
Now we show that $p$ and $H$ do not depend on the choice of the point $(r,s) \in Z$. Indeed, as $H(at) = H(bt) =
H(t)$ for every $t \in \mathbb{R} \setminus \{0\}$, it follows from (4) that $a^p + b^p =1$. This and the
relation $a,b \in (0,1)$ force that $p$, and then consequently $H$, in the representation of $f$ are unique.

We see now that $\varphi$ is the characteristic function of a semi-stable distribution (see eg. \cite{Sato},
Chapter 3) with the characteristic exponent $p$, and consequently $p \in (0,2]$.

Now take any $(r_0,s_0) \in (0,\infty)^2 \setminus Z$. Then we have
$$
\varphi(a_0 t) \varphi(b_0 t) = \varphi (c_0 t ) \psi(t), \hspace{5mm} t\in \mathbb{R},
$$
with $a_0 = {{r_0}/{d(r_0,s_0)}}$, $b_0 = {{s_0}/{d(r_0,s_0)}}$, $c_0 = {{c(r_0,s_0)}/{d(r_0,s_0)}}$. Since
$\varphi$ is positive, so is $\psi$. Define $g\colon \mathbb{R} \rightarrow \mathbb{R}$ by $ g = -\ln \psi$. The
last condition implies now that
$$
g(t) = |t|^{p} \left( a_0^{p} H(a_0 t) + b_0^{p} H(b_0 t) - c_0^{p} H(c_0 t) \right), \hspace{5mm} t\in
\mathbb{R}\setminus \{0\},
$$
so it is enough to define $K$ by $K(t) = a_0^{p} H(a_0 t) + b_0^{p} H(b_0 t) - c_0^{p} H(c_0 t)$. Observe that
$g$ is positive in a vicinity of zero. Thus, as it satisfies equation $(1)$ like $f$, it follows from Lemma 3
that $g$ is positive. Since $g(t) = |t|^p K(t)$, $t \in \mathbb{R}\setminus \{0\}$, also $K$ is positive. To
obtain the final assertion it is sufficient to insert the forms of $\varphi$ and $\psi$ into the equation. \qed

\vspace{2mm}

\noindent
\begin{cor}
Let $c,d \colon [0,\infty)^2 \rightarrow [0,\infty)$. Assume that $d(r,s) > 0$ for some $r,s >0$ and the group
generated by the set $\{ {r/s}\colon (r,s)\in Z \}$ is dense in $(0,\infty)$. Let $\varphi, \psi \colon
\mathbb{R} \rightarrow \mathbb{R}$ be nontrivial characteristic functions satisfying equation $(2)$. Then there
are $p \in (0,2]$ and $A,B >0$ such that
      $$
      \varphi(t) = \exp\{ - A|t|^{p}\}, \hspace{4mm} \psi(t) = \exp\{ - B|t|^{p}\}, \hspace{4mm} t
      \in\mathbb{R},
      $$
      and
      $$
      A \Bigl(r^{p} + s^{\alpha} - c(r,s)^{p}\Bigr) = B  d(r,s)^{p}, \hspace{5mm} r,s \geqslant 0.
      $$
\end{cor}

\vspace{2mm}

\noindent {\bf Proof.} Observe that given a function $G\colon \mathbb{R}\smallsetminus \{0\} \rightarrow
\mathbb{R}$ the set
$$
\left\{ {{r}/{s}} \colon\ r,s \in (0,\infty),\ G(rt) = G(st) \hbox{ for each } t \in \mathbb{R}\smallsetminus
\{0\} \right\}
$$
is a subgroup of the group $((0,\infty), \cdot)$. By Theorem 4 we know that $H(rt) = H(st)$ and $K(rt) = K(st)$
for all $(r,s) \in Z$ and $t\in \mathbb{R} \setminus \{0\}$. Thus the density assumption and the continuity of
the functions $H$ and $K$ imply that $H$ and $K$ are constant and inserting the form of $\varphi$ and $\psi$
into (2) implies the assertion. \qed

\vspace{2mm}

 \noindent {\large\bf Remark 3.} As follows from Remark 1 the density assumption in Corollary 1 is satisfied in each of the following
 cases concerning the set $Z$ of all zeros of the function $d$:
 \begin{namelist}{ll}
\item[-] interior of $Z$ is non-empty; %
\item[-] there are $(r_1,s_1), (r_2,s_2) \in Z$ such that ${{\ln{(r_1/s_1)}}/{\ln{(r_2/s_2)}}}$ is irrational; %
\item[-] $\lim_{n\rightarrow\infty} {{r_n}/{s_n}} =1$ for some points $(r_n,s_n) \in Z$ with $r_n \neq s_n$, $n
         \in \mathbb{N}$.
 \end{namelist}

\vspace{2mm}

\noindent {\large\bf Remark 5.} Notice that in the case of constant functions $H\equiv A$ and $K\equiv B$ we can
arbitrarily choose the function $c\colon [0,\infty)^2 \rightarrow [0,\infty)$ such that $c(r,s) \leqslant \|
(r,s)\|_{p}$, and define $d$ by
$$
d(r,s)^{p} = \frac{C}{A} \Bigl(r^{p} + s^{p} - c(r,s)^{p}\Bigr).
$$

\vspace{2mm}

\noindent {\large\bf Remark 6.} Recall that a random variable $X$ with the characteristic function $\phi$ is
semi-stable if there exist constants $r, b, c>0$ such that for each $t \in \mathbb{R}$
$$
\phi(t)^{r} = \phi\left( c t\right) e^{ibt}.
$$
Notice that if $X$ is semi-stable with constants $r,b,c$ then $X + \frac{b}{r}(c-1)$ is semi-stable with
constants $r,0,c$, thus without lost of generality, we can assume that $X$ is semi-stable if this condition
holds for $b=0$. More about semi-stable variables, including canonical form of the characteristic function,
L\'evy-Khintchine representation and examples, one can find in \cite{Sato}.

Sometimes (see eg. \cite{Bingham}) authors assume that $r$ is a natural number greater than $1$. Such definition
do not cover all the possible semi-stable variables, but it has a very natural interpretation: there exists
$r\in \mathbb{N}$ and a positive constant $c>0$ such that
$$
X_1 + \dots + X_r \stackrel{d}{=} cX,
$$
where $X_1, \dots X_r$ are independent copies of $X$. It was shown in \cite{Bingham} that branching processes
provides a large class of semi-stable distributions occurring naturally.

\subsection{$(c,d)$-pseudostable measures with respect to the weak generalized convolution}

Equation (2) was discussed by K. Oleszkiewicz in \cite{KO} in the special case when $\mu$ is a symmetric
Gaussian distribution on $\Rr$ and by J. Misiewicz and G. Mazurkiewicz (see \cite{MM}) in the case when $\mu$ is
a symmetric $p$-stable distribution on $\Rr$ for some $p\in (0,2)$. In both cases they have obtained that the
function $\varphi$ has to be of the form
$$
\varphi(t) = \varphi_{{\small C,D}, q,p}(t)= \exp\left\{ - C|t|^q - D|t|^p \right\}
$$
for some nonnegative numbers $C,D$ and some $q>0$. The problem is that not all configurations of these three
parameters are available in the sense that for $q>2$ and some numbers $C,D$ the function $\varphi_{{\small
C,D},q,p}$ is not positive definite and, consequently, it cannot be a characteristic function. Otherwise we would take $D_n \rightarrow 0$ for $n\rightarrow \infty$ and
then $\varphi_{{\small C,D}_n,q,p}(t) \rightarrow \exp\{ - C|t|^q \}$ would be positive definite, which,
however, is impossible for $q>2$. This argument shows that for any fixed $q>2$ and $C=1$ the numbers $D$ have to
be greater than some positive  number.

Maybe it is more interesting that for some parameters $C,D >0$ and some $q>2$ the function $\varphi_{{\small
C,D},q,p}$ is positive definite. It was shown in \cite{KO} and \cite{MM} that for every $r,s \geqslant 0$
$$
c(r,s) = \|(r,s)\|_q = \left(r^q + s^q \right)^{1/q}, \hspace{6mm} d(r,s)^{p} = \frac{D}{C} \left(
\|(r,s)\|_p^p - \|(r,s)\|_q^p \right).
$$
Since $d$ is  nonnegative we have to require that  $\|(r,s)\|_p \geqslant \|(r,s)\|_q$ for all $r,s
>0$, which implies that $q\geqslant p$, thus we have the following

\begin{prop}
Let $p\in(0,2]$ and $q>0$. If the pair $(\varphi_{{\small C,D},q,p}, \widehat{\gamma_p})\in \Phi_2$ is a
solution of equation $(2)$ for some $C,D >0$ then $q \geqslant p$.
\end{prop}

Nevertheless, for all $0 < p,q \leqslant 2$ and every choice of the parameters $C,D\geqslant 0$ the function
$\varphi_{{\small C,D},q,p}$ is a characteristic function as a product of a symmetric $p$-stable and a symmetric
$q$-stable characteristic functions.

\vspace{2mm}

In \cite{KO} Oleszkiewicz proved that
\begin{itemize}
 \item[$\bullet$] for $q \in (2,4] \cup \bigcup_{k=2}^{\infty}
    [4k-2, 4k]$ none of the functions $\varphi_{{\small C,D},q,2}$ can be a
    characteristic function;
 \item[$\bullet$] for $q\in \bigcup_{k=1}^{\infty}(4k,4k + 2)$ it is
    possible to find $c,d>0$ such that $\varphi_{{\small C,D},q,2}$ is a
    characteristic function.
\end{itemize}
In \cite{MM} Mazurkiewicz and Misiewicz proved that for $0< p \leqslant 1$, $q>2$ one can find $C,D>0$ such that
$\varphi_{{\small C,D},q,p}$ is a characteristic function.

Notice that in the situation considered here the information that for fixed parameters the function
$\varphi_{{\small C,D},q,p}$ is a characteristic function does not mean yet that the pair $(\varphi_{{\small
C,D},q,p}, \widehat{\gamma_p})$ is a solution of  equation $(2)$, since we have also the  assumption
$(\varphi_{{\small C,D},q,p}, \widehat{\gamma_p}) \in \Phi_2$, that is for some $\lambda \in \mathcal{P}_{+}$
$$
\varphi_{{\small C,D},q,p}(t) = \int_{[0,\infty)} \exp\left\{ - |t|^p s^p \right\} \lambda(ds), \hspace{5mm} t \in \mathbb{R}.
$$
Substituting $|t|^p \rightarrow t$  and $\alpha = {q/p}$
we would have
$$
\exp\left\{ - C t^{\alpha} - D t \right\} = \int_{[0,\infty)} \exp\left\{ - t \, s^p \right\} \lambda(ds),
\hspace{5mm} t \geqslant 0,
$$
which means that the function $[0,\infty) \ni t \rightarrow \exp\left\{ - C t^{\alpha} - D t \right\}$ would be
completely monotonic as the Laplace transform of the random variable $\theta^{p}$, where $\theta$ has the
$\gamma_p$-weakly stable distribution $\lambda$ with symmetric $p$-stable  $\gamma_p$.

\begin{prop}
Let $C,D>0$, $p\in(0,2]$ and $q\in \bigcup_{k=1}^{\infty}((2k-1)p,2kp)$. Then none of the pairs
$(\varphi_{{\small C,D},q,p}, \widehat{\gamma_p})$ is a solution of equation $(2)$.
\end{prop}

\noindent {\large\bf Proof.} Let $w(t) = C t^{\alpha} + D t$ for every $t>0$, where $\alpha = {q/p}$.  If
$(\varphi_{{\small C,D},q,p}, \widehat{\gamma_p})$ is a solution of equation $(2)$ in the class $\Phi_2$ then
the formula $g(t) = \exp\{ - w(t)\}$ defines a completely monotonic function on the positive half-line.  We see
that $g'(t) = - w'(t) g(t) < 0$ for every $t>0$. The second derivative takes the form
$$
g''(t) = g(t) \left[ (w'(t))^2 - w''(t) \right].
$$
Since $\alpha = {q/p}>1$, $w''(t) = C \alpha (\alpha -1) t^{\alpha -2}$ for every $t>0$ and $w'(0^{+})= d$ then
for $\alpha \in (1,2)$ we obtain
$$
\lim_{t\rightarrow 0^{+}} g''(t) = - \infty,
$$
which means that $g$ is not a completely monotone function.

Consider now $k\in \Nn$, $k>1$, and let $\alpha \in ((2k-1), 2k)$. The $(2k)$-derivative of the function $g$
can be written in the following form:
$$
g^{(2k)}(t) = \left[ h_{2k}(t) - w^{(2k)} \right] g(t),
$$
where $h_1(t) \equiv 0$, and for every $m \in \mathbb{N}$
$$
h_{m+1}(t) = h_m'(t) - w'(t) \left[ h_m(t) - w^{(m)}(t) \right], \hspace{5mm} w^{(m)}(t) = C \prod_{j=0}^{m-1}
(\alpha - j) \cdot t^{\alpha -m}.
$$
It is easy to see that $h_{2k}(t)$ is a linear combination of elements of the form $t^{\beta_j}$ with either
$\beta_j \geqslant \alpha -2k +1$, or $\beta_j = 0$ for every $j$, thus $h_{2k}(0^{+})$ is finite for $\alpha
\in ((2k-1), 2k)$. We see also that $w^{(2k)}(0^{+}) = + \infty$, so $g^{(2k)}(0^{+}) = - \infty$, which means
that $g$ cannot be completely monotonic, contrary to our assumptions. \qed

\vspace{3mm}

\begin{prop}
Let $C,D >0$, $p\in (0,2]$ and $q>p$. If $(\varphi_{{\small C,D},q,p}, \widehat{\gamma_p}) \in \Phi_2$ is a
solution of equation $(2)$, then for each $\alpha \in (0,1)$ the pair $(\varphi_{{\small C,D},q\alpha,p\alpha},
\widehat{ \gamma_{p\alpha}})$ is also a solution of equation {\rm (2)} and belongs to $\Phi_2$.
\end{prop}

\noindent {\large\bf Proof.} It is easy to check that the pair $(\varphi_{{\small C,D},q\alpha,p\alpha},
\widehat{ \gamma_{p\alpha}})$ satisfies equation (2). Thus we shall only prove that the function
$\varphi_{{\small C,D},q\alpha,p\alpha}$ is the mixture of the function $\widehat{ \gamma_{p\alpha}}$ with
respect to some measure from $\mathcal{P}_{+}$. To see this, observe first that by our assumptions there exists a
measure $\lambda \in \mathcal{P}_{+}$ such that
$$
\exp\left\{ - C|t|^q - D|t|^p \right\} = \int_{[0,{\infty})} \exp\left\{ - |t|^p s^p \right\} \lambda(ds),
\hspace{5mm} t \in \mathbb{R}.
$$
Let $X_p$ be a symmetric stable random variable with the characteristic function $\mathbb{R} \ni t \mapsto
\exp\{ - |t|^p\}$, $Q$ a random variable with distribution $\lambda$ and $\theta_{\alpha}$ a positive random
variable with Laplace transform $[0,\infty) \ni t \rightarrow \exp\{ - t^{\alpha} \}$ such that $X_p, Q$ and
$\theta_{\alpha}$ are independent. Then we have
\begin{eqnarray*}
\lefteqn{ \mathbb{E} \exp\left\{ i t X_p \theta_{\alpha}^{1/p}
      Q^{1/{\alpha}} \right\} = \mathbb{E} \exp\left\{ - |t|^p
      \theta_{\alpha} Q^{p/{\alpha}} \right\}} \\
 & = & \int_{[0, {\infty})} \exp\left\{ - |t|^{p\alpha} s^p \right\}
      \lambda(ds) \\
 & = & \exp\left\{ - C|t|^{q \alpha} - D|t|^{p\alpha} \right\}
      = \varphi_{{\small C,D},q\alpha,p\alpha}(t)
\end{eqnarray*}
for every $t \in \mathbb{R}$. This shows that $\varphi_{{\small C,D},q\alpha,p\alpha}$ is a scale
mixture of the function $\widehat{ \gamma_{p\alpha}}$ with respect to the distribution of the random variable
$Q^{1/{\alpha}}$, as required. \qed

\vspace{2mm}

\noindent {\large\bf Acknowledgement.} The authors are indebted to the referee for his/her valuable remarks and
comments.


\begin{thebibliography}{99}
{\small
\bibitem{BJ} Baron, K., and Jarczyk, W. (1987). {\it On a way of
             division of segments}, Aequationes Math. { 34},
             195--205.
\bibitem{Bingham} Bingham, N. H. (1988).
             {\it On the limit of a supercritical branching process, A Celebration of Applied Probability}, J.
             Appl. Probab. Special Vol. 25A, 215--228.
\bibitem{MBO2} Borowiecka--Olszewska, M. (2005). {\it The functional
              equation and strictly substable random vectors},
              Probab. Math. Statist. {25}(2),
              267--278.
\bibitem{Copson} Copson, E. T. (1935).
              {\em An Introduction to the Theory of Functions of a Complex Variable}, Oxford University Press,
              London.
\bibitem{Dud} Dudley, R.M. (1989).
              {\it Real Analysis and Probability},
              Wadsworth{\&}Brooks/Cole Advanced
              Books{\&}Software, Wadsworth Inc., Belmont,
              California.
\bibitem{FKotzNg} Fang, Kai Tai, Kotz, S. and Ng, Kai Wang. (1990).
              {\it Symmetric Multivariate and Related Distributions}, {\em Monographs on Statistics and Applied
              Probability} {36}, Chapman \& Hall, London.
\bibitem{3s} Feller, W. (1966). {\it An Introduction to Probability
             Theory and its Applications}, Vol. 2, John Wiley,
             New York.
\bibitem{Hardy Wright} Hardy, G. H. and Wright, E.M. (1979).
             {\em An Introduction to the Theory of Numbers. Fifth edition}, The Clarendon Press, Oxford
             University Press, New York.
\bibitem{Jar}Jarczyk, W. (1991). {\it A recurrent method of solving
             iterative functional equations}, Prace Naukowe
             Uniwersytetu \'{S}l\c{a}skiego w Katowicach, 1206,
             Uniwersytet \'{S}l\c{a}ski, Katowice.
\bibitem{Kingman} Kingman, J. F. C. (1963).
             {\it Random walks with spherical symmetry}, Acta Math. 109, 11--53.
\bibitem{KU} Kucharczak, J., and Urbanik, K. (1974). {\it Quasi-stable
             functions}, {Bull. Pol. Acad. Sci. Math.} {22}(3), 263-268.
\bibitem{KU2}  Kucharczak, J., and Urbanik, K. (1986). {\it
             Transformations preserving weak stability},
             Bull. Pol. Acad. Sci. Math. {34} (7-8), 475-486.
\bibitem{Ku} Kuczma, M. (1985). {\em An Introduction to the Theory of
             Functional Equations and Inequalities. Cauchy's
             Equation and Jensen's Inequalities}, Pa\'{n}stwowe Wydawnictwo
             Naukowe and Uniwersyt {\'S}l{\c{a}}ski, Warszawa
             - Krak\'{o}w - Katowice.
\bibitem{lacz} Laczkovich, M. (1986). {\it Nonnegative measurable
             solutions of difference equations}, J. London
             Math. Soc. (2) {34},   139--147.
\bibitem{4s} Luk\'{a}cs, E. (1960). {\it Characteristic Functions},
             Griffin, London.
\bibitem{MM} Misiewicz, J.K., and Mazurkiewicz, G. (2005) {\it On
             $(c,p)$-pseudostable random variables}. {J. Theoret.
             Probab.} {18(4)},   837--852.
\bibitem{MOU} Misiewicz, J.K., Oleszkiewicz, K., and Urbanik, K. (2005).
            {\it Classes of measures closed under mixing and
            convolution. Weak stability},
            {Studia Math.} 167 (3), 195--213.
\bibitem{Mis06} Misiewicz, J.K. (2006). {\it Weak stability and
            generalized weak convolution for random vectors and
            stochastic processes.}
            {\em Dynamics $\&$ Stochastics, IMS Lecture Notes
            Monograph Series} {48},   109--118.
\bibitem{KO} Oleszkiewicz, K. (2003). {\it On $p$-pseudo-stable
             random variables, Rosenthal Spaces and $\ell_p^n$
             ball slicing}, in {\em Milman V.D. and Schechtman
             G.} (eds.), {\em Lecture Notes in Math.} 1807,
             {\em Geometric Aspects of Functional Analysis}, Israel
             Seminar 2001-2002, Springer-Verlag, Berlin
             Heidelberg,   188--210.
\bibitem{JR} {Rosi{\'n}ski}, J. {\it Tempering stable
             processes}, preprint.
\bibitem{ST} Samorodnitsky, G., and Taqqu, M.S. (1994). {\it Stable
             non-Gaussian Random Processes. Stochastic Models
             with Infinite Variance.} {Chapman \& Hall}.
\bibitem{Sato} Sato, Ken-Iti. (1999). {\it L\'evy Processes and Infinitely
              Divisible Distributions}, {\em Cambridge Studies in advanced mathematics 68}, {Cambridge University
              Press}.
\bibitem{Schoenberg1}  Schoenberg, I.J. (1938). {\it Metric spaces
             and completely monotonic functions}. {Ann. of
             Math.} {38},   811--841.
\bibitem{Urbanik1}  Urbanik, K. (1976). {\it Remarks on ${\cal
             B}$-stable probability distributions}, {Bull. Pol. Acad. Sci.
             Math.} {24}(9),    783-787.
\bibitem{Urbanik64} Urbanik, K. (1964). {\it Generalized convolutions},
             Studia Math. {23},   217--245.
\bibitem{Urbanik73} Urbanik, K. (1973). {\it Generalized convolutions
             II}, Studia Math. {45},   57--70.
\bibitem{Urbanik84} Urbanik, K. (1984). {\it Generalized convolutions
             III}, Studia Math. {80},   167--189.
\bibitem{Urbanik86} Urbanik, K. (1986). {\it Generalized convolutions
             IV}, Studia Math. {83},   57--95.
\bibitem{vol1}  Vol'kovich, V. (1992). {\it On symmetric stochastic
             convolutions},  J. of Theoret.
             Probab. {5}(3),    417--430.
\bibitem{vol2}  Vol'kovich, V. (1985). {\it On infinitely
             decomposable measures in algebras  with stochastic
             convolution},  {\em Stability Problems of Stochastic
             models}. {\em Proceedings of VNIICI Seminar}, M.,
               15--24 (in Russian).
\bibitem{vol3}  Vol'kovich, V. (1984).  {\it Multidimensional ${\cal
             B}$-stable distributions and  some generalized
             convolutions}. {\em  Stability Problems of Stochastic
             Models}. {\em Proceedings of VNIICI Seminar},
               40--53 (in Russian).
\bibitem{9s} Zolotarev, V.M. {\it One-dimensional stable
             distributions}, Transl. Math. Monographs 65, Amer.
             Math. Soc., Providence.
}
\end{thebibliography}
\end{document}